\documentclass{amse-new}

\numberwithin{equation}{section} 

\usepackage{setspace}

\usepackage{amsmath, amssymb, amscd, eucal}
\usepackage{amsfonts, amsthm}
\usepackage{latexsym}
\usepackage{graphicx}
\usepackage{hyperref}
\usepackage{color}

\begin{document}


\abovedisplayskip 6pt plus 2pt minus 2pt \belowdisplayskip 6pt
plus 2pt minus 2pt
\def\vsp{\vspace{1mm}}
\def\th#1{\vspace{1mm}\noindent{\bf #1}\quad}
\def\proof{\vspace{1mm}\noindent{\it Proof}\quad}
\def\no{\no number}
\newenvironment{prof}[1][Proof]{\noindent\textit{#1}\quad }
{\hfill $\Box$\vspace{0.7mm}}
\def\q{\quad} \def\qq{\qquad}
\allowdisplaybreaks[4]


\def\endpf{\quad\rule{2mm}{2mm}}

\newcommand{\II}{\|\hspace{-0.2cm}|}
\newcommand{\tr}{\mbox{tr}}
\newcommand{\Div}{\mbox{div}}
\renewcommand{\div}{\mbox{div}}
\newcommand{\Hess}{\mbox{Hess}}
\newcommand{\Ric}{\mbox{Ric}}
\newcommand{\R}{\mathbb R}

\newcommand{\be}{\begin{equation}}
\newcommand{\ee}{\end{equation}}
\newcommand{\bee}{\begin{equation*}}
\newcommand{\eee}{\end{equation*}}
\newcommand{\bal}{\begin{aligned}}
\newcommand{\eal}{\end{aligned}}

\def\vh{\vspace{.2cm}}

\def\p{\partial}
\def\M{\mathcal{M}}
\def\la{\langle}
\def\ra{\rangle}
\def\lf{\left}
\def\ri{\right}
\def\Pi{\displaystyle{\mathbb{II}}}

\renewcommand\({\left(}
\renewcommand\){\right)}

\def\tS{\tilde{\Sigma}}
\def\m{\mathfrak{m}}
\def\mS{\mathbb{S}}
\def\tH{\tilde{H}}
\def\l{\lambda}
\def\e{\epsilon}
\def\Msg{\mathfrak{M}_{\Sigma, \gamma}}

\def\Ricg{\Ric}
\def\tg{\tilde{g}}
\def\tR{\tilde{R}}
\def\F{\mathcal{F}}
\def\E{\mathcal{E}}
\def\c{\mathfrak{c}}
\def\C{\mathcal{C}}
\def\Q{\mathcal{Q}}


\AuthorMark{Miao P.}        

\TitleMark{Monotonicity of harmonic functions on $3$-manifolds}     

\title{Monotonicity of harmonic functions on $3$-manifolds \\
with an asymptotically flat end    
}             

\author{Pengzi \uppercase{Miao}}             
    {Department of Mathematics, University of Miami, \\
     Coral Gables, FL 33146, USA \\
     E-mail\,$:$ pengzim@math.miami.edu}

\maketitle%

\Abstract{We derive monotone properties of positive harmonic functions on three dimensional manifolds 
with nonnegative scalar curvature, with an asymptotically flat end. Rigidity characterization of 
spatial Schwarzschild manifolds with two ends is also given.}


\MRSubClass{83C99, 53C20}      

\section{Introduction} 
We begin with an example that may illustrate the first two theorems in this paper. 
Given a number $ m \ne 0 $, let $ r > 0 $ be a constant satisfying 
$$  1 + \frac{m}{ 2 r }  > 0 . $$
Consider a portion of the spatial Schwarzschild manifold
$$ \left( M_r , g_m \right) =  \left( \R^3 \setminus \left\{ | x | < r  \right\} ,  \left( 1 + \frac{m}{ 2 |x| } \right)^4 g_0 \right) . $$
Here the spatial Schwarzschild metric $g_m$ is written in 
isotropic coordinates and $g_0$ denotes the Euclidean metric on $ \R^3$.

As $ m \ne 0 $, the function 
$$
u =  \frac{1}{ 1 - \left( 1 + \frac{m}{ 2 r  } \right)^{-1} } \, 
\left[ \left( 1 + \frac{m}{ 2 |x|} \right)^{-1} - \left( 1 + \frac{m}{ 2 r} \right)^{-1} \right] 
$$
is harmonic on $(M_r, g_m)$, with $ u = 0$ at  $ \Sigma =  \p M_r$ and $u \to 1$ as $ |x| \to \infty$. 

The parameter $m$ is often known as the mass of $(M_r, g_m)$. The capacity $\c_{_\Sigma} $ of the boundary 
$\Sigma$ in $(M_r, g_m)$ (see the definition in \eqref{eq-def-capacity-h}) is given by 
\be  \label{eq-c-m-S}
 \c_{_\Sigma}  = r + \frac{m}{2} .  
\ee
The ratio $ m \c_{_\Sigma}^{-1}$ satisfies
\be \label{eq-gm-ratio}
m \c_{_\Sigma}^{-1} - 2  
= -2 \left( 1 + \frac{m}{2r} \right)^{-1}. 
\ee
Given any $ t \in [0,1)$, let $ \Sigma_t = \{ u = t \} $, and 
  $$L (t) = \frac{1}{4\pi}  \int_{\Sigma_t} | \nabla u |^2  \, d \sigma , $$
where ``$ \nabla $, $ | \cdot | $, $ d \sigma$" denote the gradient, the length, the surface measure with respect to $g_m$, respectively. 
Calculation (see \eqref{eq-L-t-u-gauge} for instance) gives 
\be \label{eq-L-t-S}
L(t) =  \left( \frac{   1 + \frac{m}{2r}  t } {  1 + \frac{m}{2 r} }  \right)^2
( 1 - t )^2 .
\ee 
Therefore, if 
\be \label{eq-S-t-S}
\mathcal{S}(t) =   2 \left[ \frac{t}{1-t} - \frac{ L(t)^\frac12} { ( 1 - t)^2} \right] , 
\ee
then 
$$ \mathcal{S}(t) = - 2 \left( 1 + \frac{m}{2r}  \right)^{-1}, \ \forall \, t . $$
Together with \eqref{eq-gm-ratio}, this shows 
\be \label{eq-S-t-m-c-S}
\mathcal{S}(t) = m \c_{_\Sigma}^{-1} - 2 , \ \forall \, t . 
\ee

One may also consider the case of $ m = 0 $ 
in which $(M_r, g_0)$ is the Euclidean space $ \R^3$ minus a round ball of radius $r$.
In this case, if $u$ is the harmonic function tending to $1$ at $\infty$ and vanishing at the boundary, 
then 
$$ L(t) = (1-t)^2 , \ \forall \, t , $$
which shows the relation \eqref{eq-S-t-m-c-S} 
holds on $(M_r, g_0)$ as well. 

It turns out that \eqref{eq-S-t-m-c-S} 
reveals general properties 
of positive harmonic functions on a class of $3$-manifolds with nonnegative scalar curvature.  
To state our results, we recall a few definitions first. 

A Riemannian $3$-manifold $(M, g)$ is asymptotically flat (AF)  if $M$,  
outside a compact set, is diffeomorphic to $ \R^3$ minus a ball;
the associated metric coefficients  satisfy 
$$
g_{ij} = \delta_{ij} + O ( |x |^{-\tau} ), \ 
 \p g_{ij} = O ( |x|^{-\tau -1}) ,  \ \ \p \p g_{ij} = O (|x|^{-\tau -2} ), 
$$
for some $ \tau > \frac12$; and the scalar curvature  of $g$ is  integrable.
Under these AF conditions, the limit 
$$
\m  = \lim_{ r \to \infty  } \frac{1}{16 \pi} \int_{ |x | = r } \sum_{j, k} ( g_{jk,j} -  g_{jj,k} ) \frac{x^k}{ |x| } 
$$
exists and is called the mass of $(M,g)$ (see  \cite{ADM61, SchoenYau79, Witten81} for instance). 
It is known (\cite{Bartnik86, Chrusciel86}) 
that $\m$ is a geometric quantity, independent on the choice of the coordinates $\{ x_i \}$.

On an asymptotically flat $3$-manifold $(M,g)$ with boundary $ \Sigma  $, 
the capacity (or $L^2$-capacity) of $\Sigma $ is defined by 
$$
\c_{_\Sigma} = \inf \left\{  \frac{1}{4 \pi} \int_M | \nabla f |^2  \ : \ 
\, f |_{\Sigma} = 1 \ \text{and} \  f \to 0 \ \text{at} \ \infty \right\} .
$$
Equivalently, if $\phi$ denotes the function with 
$$
\Delta \phi = 0 , \ \ 
 \phi |_\Sigma = 1 , \ \text{and} \  
 \phi \to 0  \ \ \text{at} \  \ \infty , 
$$
then $ \c_{_\Sigma} $ appears as the coefficient of $|x|^{-1}$ in the expansion 
\be \label{eq-def-capacity-h}
 \phi = \c_{_\Sigma} |x|^{-1} + o ( | x|^{-1} ) , \ \  \text{as} \ |x| \to \infty . 
\ee

\begin{theorem} \label{thm-main-ineq}
Let $(M, g)$ be a complete, orientable, asymptotically flat $3$-manifold
with nonnegative scalar curvature,  with boundary $\Sigma$.
Suppose $ \Sigma$ is connected and $H_2 (M, \Sigma) = 0$.
Let $ u $ be the harmonic function on $(M,g)$ with $ u = 0 $ at $ \Sigma$ and $ u \to 1$ at $\infty$. 
For each regular value $ t \in [0,1)$ of $u$, define $  \mathcal{S}(t)   $ by 

\be \label{eq-def-St-main}
 \displaystyle   \mathcal{S}(t)  
= \left\{
\begin{split}
2 \left[ \frac{t}{(1-t)} - \frac{L(t)^\frac12} { (1-t)^2} \right], & \ \ \text{if} \  L(t)  > t^2 (1-t)^2  , \\
\vh & \\ 
\frac{t}{(1-t) } \left[ 1 -  \frac{L(t)}{ t^2 ( 1 - t)^2 } \right], & \  \ \text{if} \  L(t)   \le  t^2 (1-t)^2  . 
\end{split}
\right.
\ee 
Here $ \displaystyle  L (t) = \frac{1}{4\pi} \int_{ u^{-1} (t) } | \nabla u |^2 \, d \sigma $,
where ``$ \nabla $, $ | \cdot | $, $ d \sigma$" denote the gradient, the length, the surface measure with respect to 
the metric $g$, respectively.

Then $\mathcal{S}(t)$ is monotone nondecreasing in $t$, i.e.
$$
   \mathcal{S} (t)   \nearrow  \ \text{as} \ t \nearrow , 
$$
and 
$$
 \mathcal{S}(t) \le \m \c_{_\Sigma}^{-1} - 2 , \ \forall \, t . 
$$
Here $ \m $  is the mass of $ (M, g)$, and $\c_{_\Sigma} $ is the capacity of $\Sigma$, which is 
the coefficient in the expansion 
$$
u = 1 -  \frac{\c_{_\Sigma}}{ | x| } + o ( | x|^{-1} ) , \ \text{as} \ |x| \to \infty. 
$$
Moreover, 
if $  \mathcal{S}(t) = \m \c_{_\Sigma}^{-1} - 2 $ at some $ t \in [0,1) $,   
then the region $\{ u \ge t \}$ in $(M,g)$ is isometric to a spatial Schwarzschild manifold outside a rotationally symmetric sphere.  
\end{theorem}

\vspace{.1cm}

\begin{remark}
If  $  (M, g) = (M_r, g_m)$, then $ L(t) > t^2 ( 1-t)^2$ by \eqref{eq-L-t-S}. 
\eqref{eq-def-St-main} explains the choice of $\mathcal{S} (t) $ in  \eqref{eq-S-t-S}. 
\end{remark}

\begin{remark}
In Theorem \ref{thm-main-ineq}, at $ t = 0 $,  $L(0) > 0 $ by the Hopf lemma. Hence, 
\be \label{eq-S-t-small-t}
   \mathcal{S}(t)  = 2 \left[ \frac{t}{(1-t)} - \frac{L(t)^\frac12} { (1-t)^2} \right] 
\ \ \text{for small } t > 0 . 
\ee
As $ t \to 1 $, if $ \p \p \p g_{ij} = O (|x|^{-\tau - 3} ) $ at $\infty$, then by \cite[Theorem 2.1 (ii)]{M22},  
\be \label{eq-limit-L-q}
\lim_{t \to 1} \frac{1}{1-t} \left[ 1 - \frac{L(t) }{(1-t)^2} \right] = \m \c^{-1}_{_\Sigma} . 
\ee
This shows that a comparison between 
$ L (t)$ and $ t^2 ( 1 - t)^2 $ near $\infty$ 
is suggested by comparing $ \m \c_{_\Sigma}^{-1}$ and $ 2$.
More precisely, by \eqref{eq-limit-L-q}, 
$$ \frac{1}{1-t} \left[ \frac{L(t)}{(1-t)^2} - t^2 \right] = 1 + t - \m \c_{_\Sigma}^{-1} + o(1) . $$
Consequently, 
$$ \text{if} \ \m \c_{_\Sigma}^{-1} > 2 , \ \text{then} \ 
 L (t) < t^2 ( 1 - t)^2  \ 
\text{for} \  t  \ \text{close to} \  1 . 
$$ 
Hence, on manifolds $(M,g)$ satisfying $ \m \c_{_\Sigma}^{-1} > 2 $, 
\be \label{eq-S-t-t-close-to-1}
 \mathcal{S}(t)  = \frac{t}{(1-t) } \left[ 1 -  \frac{L(t)}{ t^2 ( 1 - t)^2 } \right] \ 
\text{for} \  t  \ \text{close to} \  1 . 
\ee
\end{remark}

\vspace{.1cm} 

\begin{theorem} \label{thm-main-2}
Let $(M, g)$ be a complete, orientable, asymptotically flat $3$-manifold
with nonnegative scalar curvature,  with boundary $\Sigma$.
Suppose $ \Sigma$ is connected and $H_2 (M, \Sigma) = 0$.
Let $ u $ be the harmonic function on $(M,g)$ with $ u = 0 $ at $ \Sigma$ and $ u \to 1$ at $\infty$. 
Suppose the mass $ \m$ of $(M,g)$ and the capacity $ \c_{_\Sigma}$ of $ \Sigma$ in $(M,g)$ 
satisfy 
$$ \m \c_{_\Sigma}^{-1}  < 2 . $$
For each regular value $ t \in [0,1)$  of $u$, let
 $  \displaystyle L (t) = \frac{1}{4\pi} \int_{ u^{-1} (t) } | \nabla u |^2 \, d \sigma $.
Then
\begin{enumerate}
\item[(i)] $ L (t) > t^2 ( 1 - t)^2 $, $ \forall \, t $; 
\item[(ii)] the function $ \mathcal{S} (t)$, defined by 
$$ \mathcal{S}(t): = 2 \left[ \frac{t}{1-t} - \frac{L(t)^\frac12}{ ( 1 - t)^2} \right] ,  $$
is monotone nondecreasing in $t$, i.e.
$$
   \mathcal{S} (t)   \nearrow  \ \text{as} \ t \nearrow , 
$$
and 
$$
 \mathcal{S}(t) \le \m \c_{_\Sigma}^{-1} - 2 , \ \forall \, t. 
$$
Moreover, 
if $  \mathcal{S}(t) = \m \c_{_\Sigma}^{-1} - 2 $ at some $ t \in [0,1) $, then
the region $\{ u \ge t \}$ in $(M,g)$ is isometric to a spatial Schwarzschild manifold outside a rotationally symmetric sphere. 
\end{enumerate} 
\end{theorem}

The benchmark $ t^2 ( 1 - t)^2$ arises from the corresponding quantity on 
a spatial Schwarzschild manifold with two ends. Given a constant $ m > 0 $, the complete spatial Schwarzschild $3$-manifold 
$(M_m, g_m)$ takes the form of  
$$ (M_m, g_m)   = \left( \R^3 \setminus \left\{ 0 \right\},   \left( 1 + \frac{m}{ 2 |x|} \right)^4 g_0  \right). $$
On $(M_m, g_m)$,  the function 
$  
u =  \left( 1 + \frac{m}{ 2 |x|} \right)^{-1}
$
is harmonic with $u \to 1$ as $ | x | \to \infty$ and $ u \to 0$ as $ |x| \to 0$. 
For each $ t \in (0,1)$, if 
$  L (t) = \frac{1}{4\pi}  \int_{ u^{-1} (t) } | \nabla u |^2  \, d \sigma , $
it can be checked 
\be \label{eq-fact-Schwarzs-m}
\left\{ 
\begin{split}
 L(t) = &  \ t^2 ( 1 - t)^2   , \ \forall \, t, \\
 \ m \, = & \ 2 c .
 \end{split}
 \right. 
\ee
Here $ c > 0  $ is the coefficient of $ - |x|^{-1}$ in the expansion of $u$ as $ |x| \to \infty$. 

\vh

The next result generalizes \eqref{eq-fact-Schwarzs-m} to spaces  that are 
analogous to $ (M_m, g_m)$. For the sake of simplicity,  
we assume the underlying manifold is $ \R^{3} \setminus \{ 0 \}$.  

\begin{theorem} \label{thm-main-app}
Let $ g $ be a Riemannian metric with nonnegative scalar curvature on $M = \R^{3} \setminus \{ 0 \} $.
Suppose
\begin{itemize}
\item 
$(M,g)$ is  asymptotically flat as $ | x | \to \infty$;  
\item 
 $ (M,g)$ admits a harmonic function $ u$ with 
$ u \to 1 $ as $ |x|  \to \infty$ and $ u \to 0 $ as $ | x | \to 0$. 
\end{itemize} 
Let $ \m $ be the mass of $(M,g)$ at the end where $  | x | \to \infty$.
Let $\c$ be
the capacity of $(M,g)$, i.e. the positive 
constant coefficient in the expansion of 
$$ u = 1 - \frac{\c}{|x|} + o ( | x|^{-1} ) , \ \text{as} \ | x |  \to \infty. $$
For each regular value $ t \in (0,1)$ of $u$, let $ \Sigma_t = u^{-1} (t)$,  
and define $ \Q(t)$ by 
$$
\Q(t)  =   \frac{t}{1-t}  \left[ 1 -  \frac{ L(t ) }{  t^2 (1-t)^2 }  \right] ,
$$
where $ \displaystyle L (t) = \frac{1}{4\pi}  \int_{\Sigma_t} | \nabla u |^2  \, d \sigma  $. 
Then the following holds.

\begin{enumerate}

\item[(i)] $ \Q(t)$ is monotone nondecreasing in $ t$, i.e. 
$ \Q(t) \, \nearrow$ as $  t \, \nearrow  $; and  
$$  \Q(t) \le \m \c^{-1} - 2 , \ \forall \,  t .  $$
Moreover,  the equality holds at some $t \in (0,1) $  if and only if $ \m = 2 \c > 0 $ and 
the region $\{ u \ge t \}$ in $(M,g)$ is isometric to  
$$ \left(   \R^3 \setminus \{ | x |  < r \} , \left(  1 +  \frac{ \m  }{ 2   | x | } \right)^4 g_0 \right) , 
\ \text{with} \   \left( 1 + \frac{m}{2 r} \right)^{-1}  = t . $$ 


\item[(ii)]  $ \displaystyle  \lim_{t \to 0+}  \frac{L(t)}{t}  $ exists (finite or $\infty$), and 
$$
 \m \c^{-1} - 2 + \lim_{t \to 0} \frac{ L(t)}{t}  \ge 0  .
$$

\vspace{.1cm} 

\noindent Moreover, the equality holds if and only if $\m > 0$ and $(M, g)$ is isometric to the spatial Schwarzschild manifold $(M_\m, g_\m)$ with two ends. 

\vspace{0.2cm}

\item[(iii)] If $  \displaystyle  \lim_{ t \to 0 + } L(t)  = 0 $ and $ \m = 2 \c $, 
then  $(M, g)$ is isometric to the spatial Schwarzschild manifold $(M_\m, g_\m)$ with two ends. 

\end{enumerate} 

\end{theorem}

\vspace{.2cm}

\begin{remark}
The metric $g$ is not assumed to be complete as $x \to 0 $ in Theorem \ref{thm-main-app}. 
For instance, $ g$ may just  be a reparametrization of a metric on $\R^3 \setminus \{ |x | <  r \}$, 
$ r > 0$, from the open interior 
$ \R^3 \setminus \{ |x | \le r \} $ to $ \R^3 \setminus \{ 0 \}$. 
\end{remark} 

\vspace{.2cm}

\begin{remark}
Under the assumption of 
 $  \displaystyle \lim_{ t \to 0 + } L(t)  = 0  $ in (iii) of Theorem \ref{thm-app-1},
it was known $ \m \ge 2 \c $ (see \cite{M24}). 
\end{remark}

\vspace{.2cm}

In recent years, there has been a sequence of works in the literature that 
explores properties of the level sets of harmonic functions under a scalar curvature condition
on a $3$-manifold. Interested readers are referred to
\cite{Stern19, BKKS19, MW21, AMO21, CL21, MW22, AMMO22, M22, O22, XYZ23}
for instance. 
These works are tied to and preceded by the study of  
the level sets of harmonic functions under a Ricci curvature condition \cite{CM13, CM14}, 
starting with Colding’s work \cite{Colding12}. 

We make use of results in \cite{M22} in the derivations of  
Theorems \ref{thm-main-ineq} -- \ref{thm-main-app}.

\section{Proof of Theorems \ref{thm-main-ineq} and \ref{thm-main-2}}

We recall the statement of Theorem \ref{thm-main-ineq}. 

\begin{theorem} \label{thm-sec-ineq-2}
Let $(M, g)$ be a complete, orientable, asymptotically flat $3$-manifold
with nonnegative scalar curvature,  with boundary $\Sigma$.
Suppose $ \Sigma$ is connected and $H_2 (M, \Sigma) = 0$.
Let $ u $ be the harmonic function on $(M,g)$ with $ u = 0 $ at $ \Sigma$ and $ u \to 1$ at $\infty$. 
For each regular value $ t \in [0,1)$ of $u$, define $  \mathcal{S}(t)   $ by 

\be \label{eq-def-St-main-2}
 \displaystyle   \mathcal{S}(t)  
= \left\{
\begin{split}
2 \left[ \frac{t}{(1-t)} - \frac{L(t)^\frac12} { (1-t)^2} \right], & \ \ \text{if} \  L(t)  > t^2 (1-t)^2  , \\
\vh & \\ 
\frac{t}{(1-t) } \left[ 1 -  \frac{L(t)}{ t^2 ( 1 - t)^2 } \right], & \  \ \text{if} \  L(t)   \le  t^2 (1-t)^2  . 
\end{split}
\right.
\ee 
Here $ \displaystyle  L (t) = \frac{1}{4\pi} \int_{ u^{-1} (t) } | \nabla u |^2 \, d \sigma $,
where ``$ \nabla $, $ | \cdot | $, $ d \sigma$" denote the gradient, the length, the surface measure with respect to 
the metric $g$, respectively.

Then $\mathcal{S}(t)$ is monotone nondecreasing in $t$, i.e.
$$
   \mathcal{S} (t)   \nearrow  \ \text{as} \ t \nearrow , 
$$
and 
$$
 \mathcal{S}(t) \le \m \c_{_\Sigma}^{-1} - 2 , \ \forall \, t. 
$$
Here $ \m $  is the mass of $ (M, g)$ and $\c_{_\Sigma} $ is the capacity of $\Sigma$ in $(M,g)$. 
Moreover, 
if $  \mathcal{S}(t) = \m \c_{_\Sigma}^{-1} - 2 $ at some $ t \in [0,1) $, 
then $\{ u \ge t \}$ in $(M,g)$ is isometric to a spatial Schwarzschild manifold outside a rotationally symmetric sphere.  
\end{theorem}

\begin{prof}
We first follow a construction in \cite[Section 7]{M22}.
Given any positive harmonic function $v$ on $(M, g)$ with $ v \to 1 $ at $\infty$, 
consider 
\begin{itemize}
\item the metric $ \bar g = v^4 g $, which has nonnegative scalar curvature, and is asymptotically flat at $  \infty$; 
\item the function $ \bar u = v^{-1} u $, which is harmonic on $(M, \bar g)$  with 
$ \bar u \to 1 $ at $ \infty $ and $ \bar u = 0 $ at $ \Sigma = \p M$. 
\end{itemize}
\noindent 
Applying \cite[Theorem 3.2 (ii)]{M22} to $(M, \bar g, \bar u)$,  we have the following monotone quantity
\be \label{eq-def-barBs}
 \mathcal{ \bar B} ( s ) = \frac{1}{1 - s } \left[ 4 \pi -   \frac{1}{ (1 - s  )^2 } \int_{ \{ \bar u = s \} }  | \bar \nabla \bar u |_{\bar g} ^2 \, d \bar \sigma \right] 
\  \  \nearrow  \ \  \text{as} \  \ s \nearrow ,
\ee
\noindent where $s  \in [0 , 1)$ is a regular value of $\bar u$, and ``$ \bar \nabla $, $ | \cdot |_{\bar g}$, $ d \bar \sigma$" denote the gradient, the length, the surface measure with respect to $\bar g$, respectively. 

\vspace{.1cm}

Given any constant $k > 0 $, we may choose 
\be \label{eq-def-v}
v = u +  k^{-1} ( 1 - u ) . 
\ee
The associated $\bar u$ is given by 
\be \label{eq-baru-and-u}
\bar u =  \frac{ 1}{ 1 +  k^{-1} (  u^{-1} - 1 ) }.
\ee
Given a regular value $ s $ of $\bar u$, 
let $ t $ be chose so that $ \{ \bar u = s \} = \{ u  = t \} $. 
By \eqref{eq-baru-and-u},
\bee
s = \frac{ t }{ t  +   k^{-1}  ( 1 - t ) }  \ \ 
\text{and} \ \ 
1 - s = \frac{ ( 1 - t )  }{ k t  + ( 1 - t ) }. 
\eee
Direct computation gives 
\be \label{eq-conformal-bdry-relation-bar}
\int_{ \{ \bar u = s \} } | \bar \nabla \bar u |_{\bar g}^2 \, d \bar \sigma = \int_{ \{ u = t \} }  | \nabla \bar u |^2 \, d \sigma ,
\ee
\be \label{eq-tw0-grads-baru}
| \nabla \bar u  | =  \left[ k u + ( 1 - u ) \right]^{-2} \, k | \nabla u | . 
\ee
It follows from \eqref{eq-def-barBs} -- \eqref{eq-tw0-grads-baru} that 
\be \label{eq-Bs-Bkt}
\mathcal{\bar B} (s) = \mathcal{B}_k (t), 
\ee
where 
\be
 \mathcal{B}_k (t)  =  \frac{ k t  + ( 1 - t ) } { ( 1 - t )  } 
\left[ 4 \pi - \frac{ k^2 }{ ( 1 - t )^2  \,  ( k t +  ( 1 - t ) )^2 }  \int_{ \{ u = t \} }  | \nabla u |^2  \, d \sigma \right] .
\ee

\vh

\noindent The monotone property of $\mathcal{\bar B} (s)$ in \eqref{eq-def-barBs}  translates to 
that of $ \mathcal{B} (t)$, that is 
\be \label{eq-monotone-Bkt}
\mathcal{B}_k (t) \   \nearrow  \  \text{as} \  t  \ \nearrow . 
\ee

\vh

By  \cite[Theorem 3.2 (ii) and Remark 3.5]{M22}, for any $ s $, 
\be \label{eq-limit-Bkt}
 \mathcal{\bar B} (s ) \le 4 \pi \m (\bar g) { \bar \c }^{-1}  ,
\ee
where $ \m (\bar g)$ is the mass of $ \bar g$, and $\bar \c $ is the coefficient in the expansion
\bee
\bar u = 1 - \bar \c | x|^{-1} + o ( | x|^{-1} ) , \ \text{as} \ |x| \to \infty. 
\eee
By \eqref{eq-def-v} and \eqref{eq-baru-and-u}, 
\be  \label{eq-mbar-cbar}
\m (\bar g) = \m - 2 ( 1 - k^{-1} ) \c_{_\Sigma} \ \ \text{and} \ \ 
\bar \c = k^{-1} \c_{_\Sigma}  . 
\ee
It follows from \eqref{eq-Bs-Bkt}, \eqref{eq-limit-Bkt} and  \eqref{eq-mbar-cbar} that
\be \label{eq-Bkt}
 \frac{1}{k} \left[  \frac{1}{4\pi}  \mathcal{B}_k (t)  - 2 \right]   \le \m \c_{_\Sigma}^{-1} - 2 .
\ee 
In terms of $ \displaystyle L(t) = \frac{1}{4\pi}  \int_{ \{ u = t \} }  | \nabla u |^2  \, d \sigma  $,
the left side of \eqref{eq-Bkt} becomes 
\be \label{eq-Bkt-2}
\begin{split}
 \frac{1}{k} \left[ \frac{1}{4\pi}   \mathcal{B}_k (t)  - 2 \right] 
= & \  \frac{  t  } {  1 - t  } -  k^{-1} - \frac{ k } {   [ k t +  ( 1 - t ) ] \,   }  \frac{ L(t)}{ ( 1 - t )^3  }. 
\end{split} 
\ee

\vh

By the monotonicity of $\mathcal{B}_k (t)$, for any $ k > 0$
and any  $t_1 , t_2 \in [0,1) $ with $ t_1 <  t_2 $, 
\be \label{eq-k-t1-t2}
 \frac{1}{k} \left[ \frac{1}{4\pi}   \mathcal{B}_k (t_1)  - 2 \right]  
\le 
\frac{1}{k} \left[ \frac{1}{4\pi}   \mathcal{B}_k (t_2)  - 2 \right] .
\ee
We now take the supremum over $k$ of both sides of \eqref{eq-k-t1-t2}. 
For each fixed $t$, define 
\be \label{eq-def-St-sup}
\mathcal{S}(t) : =    \sup_{k > 0 } \frac{1}{k} \left[ \frac{1}{4\pi}   \mathcal{B}_k (t)  - 2 \right].  
\ee
It follows from  \eqref{eq-k-t1-t2} and \eqref{eq-Bkt} that  
\be  \label{eq-St-mono-ineq}
 \mathcal{S}(t_1)  \le   \mathcal{S}(t_2)
\le  \m \c^{-1}_{_\Sigma} - 2  . 
\ee


We next show that  $\mathcal{S}(t)$ is given by \eqref{eq-def-St-main-2}. 
By \eqref{eq-Bkt-2}, 
$$ 
\mathcal{S}(t)  =   \frac{  t  } { ( 1 - t )  } - \inf_{k > 0} \, \Psi (k) .
$$
Here
$$ \Psi (k) = k^{-1} + \frac{k}{ [ k t +  ( 1 - t ) ] }  \, \frac{ L(t) } {( 1 - t )^3}     $$
and
\bee
\begin{split}
\Psi'(k) = & \ - k^{-2} +\frac{1}{  [ k t +  ( 1 - t ) ]^2 } \,  \frac{ L(t) } {( 1 - t )^2}  \\
& \ \\ 
= & \ \frac{1}{[ k t +  ( 1 - t ) ]^{2} } 
\left[ \frac{ L(t)^\frac12 } {( 1 - t )} + \frac{ k t +  ( 1 - t ) }{k} \right] 
\left[ \frac{ L(t)^\frac12 } {( 1 - t )}  -  \frac{ k t + ( 1 - t ) }{k} \right] . 
\end{split}
\eee
By elementary calculus, 
\begin{itemize}
\item if $ \displaystyle  \frac{ L(t)^\frac12 } {( 1 - t )}  > t $, then 
$ \inf_{k} \Psi(k)  $ is uniquely achieved at $ k = k_0 $, where 
\be \label{eq-k-k-0}
 \frac{1-t}{k_0} =  \frac{ L(t)^\frac12 } {1 - t } - t ; 
\ee
in this case,
$$  \inf_{k>0} \Psi(k) = \Psi (k_0) = 
\frac{1}{(1-t)} \left[  \frac{ L(t)^\frac12 } {1 - t } - t  \right] + \frac{L(t)^\frac12} { (1-t)^2} ,
$$
and consequently, 
$$ \mathcal{S} (t)  = 2 \left[ \frac{t}{(1-t)} - \frac{L(t)^\frac12} { (1-t)^2} \right]  ; $$

\item if $ \displaystyle  \frac{ L(t)^\frac12 } {( 1 - t )}  \le  t $, then $ \Psi'(k) < 0 $ for all $ k > 0 $; in this case,
$$ \inf_{k} \Psi(k) = \lim_{k \to \infty} \Psi (k) = \frac{L(t)}{ t   ( 1 - t)^3  } , $$
and 
$$ \mathcal{S} (t) =  \frac{t}{(1-t) } -  \frac{L(t)}{ t ( 1 - t)^3 }  . $$

\end{itemize}


By now, we have shown $\mathcal{S}(t)$, given by \eqref{eq-def-St-main-2}, 
satisfies \eqref{eq-St-mono-ineq}. 


Next, suppose $ \mathcal{S}(t) = \m \c_{_\Sigma}^{-1} - 2 $ 
at some $ t \in [0, 1)$. 
First consider the case where this $t$ satisfies $L (t) > t^2 ( 1 - t)^2 $.
In this case, \eqref{eq-Bkt} becomes an equality 
\bee 
 \frac{1}{k} \left[  \frac{1}{4\pi}  \mathcal{B}_k (t)  - 2 \right]  =  \m \c_{_\Sigma}^{-1} - 2 
\eee 
with $ k = k_0$ given in \eqref{eq-k-k-0}.  
Accordingly, \eqref{eq-limit-Bkt} is an equality 
$$
 \mathcal{\bar B} (s ) =  4 \pi \m (\bar g) { \bar \c }^{-1}  . 
$$ 
By \cite[Theorem 3.2 (ii)]{M22},  this occurs if and only if $(E , \bar g)$ is isometric to 
 $$ \left( \R^3 \setminus \{ | x |  < r \} , g_0 \right) $$
 for some $ r > 0$,
 where $ E$ is the region in $M$ given by 
 $$ E = \{ \bar u \ge s \}  = \{ u \ge t \}  . $$ 
In terms of the coordinate $x$, the functions $\bar u $ and $v$  on $E$ are determined  by 
\bee 
  \bar u = 1 - \frac{  ( 1 - s)  r}{ | x | }, 
\ \ \mathrm{and} \ \ 
  v^{-1} = 1 - ( 1 - k_0) \frac{( 1 - s) r}{ | x | }  .  
\eee
As a result, 
$$   g = v^{-4} g_0 = \left(  1 +   \frac{ ( k_0 - 1) ( 1 - s) r  }{  | x | } \right)^4 g_0  . $$
This shows $(E, g)$ is isometric to the part  
$$ \left(   \R^3 \setminus \{ | x |  < r \} , \left(  1 +  \frac{ \m  }{ 2   | x | } \right)^4 g_0 \right) $$
of a spatial Schwarzschild manifold of mass  
$ \m =  2  ( k_0  - 1) ( 1 - s) r     $.

We are left with the case of $ \mathcal{S}(t) = \m \c_{_\Sigma}^{-1} - 2 $ at some $t$ with $L (t) \le  t^2 ( 1 - t)^2 $.
This necessarily implies $ t > 0$. 
A proof of this case is verbatim to that of $ \Q (t) =  \m \c_{_\Sigma}^{-1} - 2  $ 
in Theorem \ref{thm-main-app} (i), which we refer to 
Section \ref{sec-main-app}. 
\end{prof} 

\begin{remark}
At $ t = 0 $, $ \displaystyle L (0 ) = \frac{1}{4\pi} \int_{\Sigma} | \nabla u |^2 \, d \sigma > 0 $. 
Hence,  
$$ \mathcal{S}(0) = - 2  L(0)^\frac12 . $$
Therefore, Theorem \ref{thm-sec-ineq-2} implies
\be \label{eq-7322}
1 - \left( \frac{1}{4\pi} \int_{\Sigma} | \nabla u |^2 \, d \sigma \right)^\frac12  \le \frac12 \m \c_{_\Sigma}^{-1} ,
\ee
which was \cite[Theorem 7.3]{M22}.
\end{remark}

\vspace{0.1cm} 

\begin{remark} \label{rem-Q}
For $ 0 < t_1 < t_2 < 1 $, 
letting $ k \to \infty$ in  \eqref{eq-k-t1-t2}, we also obtain  
\be  \label{eq-Q-bd}
 \frac{  t_1  } {  1 - t_1  }  - \frac{  L(t_1) } {  t_1 ( 1 - t_1 )^3    }  
 \le  \frac{  t_2  } {  1 - t_2  }  - \frac{  L(t_2) } {  t_2  ( 1 - t_2 )^3    }   . 
\ee
This shows, 
if $(M,g)$ and $u$ are given  in Theorem \ref{thm-sec-ineq-2}, 
the quantity
$$
Q(t)  : =   \frac{t}{1-t}  \left[ 1 -  \frac{ L(t ) }{  t^2 (1-t)^2 }  \right] , \ t > 0, 
$$
is always monotone nondecreasing in $ t$, 
and satisfies 
$
 Q(t) \le \m \c_{_\Sigma}^{-1} - 2 , \ \forall \, t .   
$
\end{remark}

\vh

Next, we give a sufficient condition that guarantees 
$ L (t) > t^2 ( 1 - t)^2 , \ \forall \, t . $

\begin{theorem} \label{thm-sec-L-t-m-c}
Let $(M,g)$, $\Sigma$ and $u$ be given as in Theorem \ref{thm-sec-ineq-2}.
Let $ \c_{_\Sigma}$ denote the capacity of $ \Sigma $ in $(M,g)$.
If $ \m $ and $ \c_{_\Sigma}  $ satisfy 
$$ \m \c_{_\Sigma}^{-1} <  2 , $$ 
then
$$ \frac{1}{4\pi} \int_{ u^{-1} (t) } | \nabla u |^2 \, d \sigma> t^2 ( 1 - t)^2 , $$ 
for any regular values $ t $.
\end{theorem}

\begin{prof} Let $ t \in [0,1)$ be a regular value of $ u$. Let $\Sigma_t = u^{-1} (t)$, 
$M_t = \{ u \ge t \}$, and 
$$ u_t = \frac{ u - t }{1 - t} . $$  
Applying \eqref{eq-7322} to 
the triple $(M_t, g, u_t)$, we have 
$$
\frac{ \m}{ 2 \c_{_{\Sigma_t} } } \ge 1 - \left( \frac{1}{4\pi} \int_{\Sigma_t} | \nabla u_t |^2 \, d \sigma \right)^\frac12 ,
$$
where  $ \c_{_{\Sigma_t} }  = ( 1 - t)^{-1} \c_{_\Sigma}  $ is the capacity of $\Sigma_t$ in $(M,g)$.

In terms of $u$, this shows
\be \label{eq-cor-of-good-result-22}
 \frac{ L(t)^\frac12  }{1-t} \ge 1 - \frac{ \m}{ 2 \c_{_{\Sigma} } }  ( 1 - t)  , 
\ee
where $ \displaystyle L (t) =  \frac{1}{4\pi} \int_{\Sigma_t} | \nabla u |^2 \, d \sigma  $.
It follows from \eqref{eq-cor-of-good-result-22} that
$$ 
 \m \c_{_\Sigma}^{-1} < 2 \ \Longrightarrow \ L(t)^\frac12 >  t (1-t)  . 
$$
\end{prof}

\begin{remark}
If $(M,g) = (M_r, g_m)$, the portion of a spatial Schwarzschild manifold outside a rotationally symmetric sphere, then $ \m \c_{_\Sigma}^{-1} < 2 $ by \eqref{eq-gm-ratio}. 
\end{remark}

\vh

Theorem \ref{thm-main-2} follows from Theorems \ref{thm-sec-ineq-2} and \ref{thm-sec-L-t-m-c}. 

\section{Proof of Theorem \ref{thm-main-app}} \label{sec-main-app}

We recall the statement of Theorem \ref{thm-main-app}.

\begin{theorem} \label{thm-app-1}
Let $ g $ be a Riemannian metric with nonnegative scalar curvature on $M = \R^{3} \setminus \{ 0 \} $.
Suppose
\begin{itemize}
\item 
$(M,g)$ is  asymptotically flat as $ | x |  \to \infty$;  
\item 
 $ (M,g)$ admits a harmonic function $ u$ with 
$ u \to 1 $ as $ | x |  \to \infty$ and $ u \to 0 $ as $ | x |  \to 0$. 
\end{itemize} 
Let $ \m $ be the mass of $(M,g)$ at the end where $ | x |  \to \infty$.
Let $\c$ be
the capacity of $(M,g)$, i.e. the positive 
constant coefficient in the expansion of 
$$ u = 1 - \frac{\c}{|x|} + o ( | x|^{-1} ) , \ \text{as} \ | x |  \to \infty. $$
For each regular value $ t \in (0,1)$ of $u$, let $ \Sigma_t = u^{-1} (t)$,  
and define $ \Q(t)$ by 
$$
\Q(t)  =   \frac{t}{1-t}  \left[ 1 -  \frac{ L(t ) }{  t^2 (1-t)^2 }  \right] ,
$$
where $ \displaystyle L (t) = \frac{1}{4\pi}  \int_{\Sigma_t} | \nabla u |^2  \, d \sigma  $. 
Then the following holds.

\vspace{.1cm} 

\begin{enumerate}

\item[(i)] $ \Q(t)$ is monotone nondecreasing in $ t$, i.e. 
$ \Q(t) \, \nearrow$ as $  t \, \nearrow  $; and  
$$  \Q(t) \le \m \c^{-1} - 2 , \ \forall \,  t .  $$
Moreover,  the equality holds at some $t \in (0,1) $  if and only if $ \m = 2 \c  > 0 $ and 
the region $\{ u \ge t \}$ in $(M,g)$ is isometric to  
$$ \left(   \R^3 \setminus \{ | x |  < r \} , \left(  1 +  \frac{ \m  }{ 2   | x | } \right)^4 g_0 \right) , 
\ \text{with} \   \left( 1 + \frac{m}{2 r} \right)^{-1}  = t . $$ 

\vspace{.1cm} 

\item[(ii)]  $ \displaystyle  \lim_{t \to 0+}  t^{-1} L(t)  $ exists (finite or $\infty$), and 
$$
 \m \c^{-1} - 2 + \lim_{t \to 0} t^{-1} L(t)  \ge 0  .
$$

\vspace{.1cm} 

\noindent Moreover, the equality holds if and only if $\m > 0$ and $(M, g)$ is isometric to the spatial Schwarzschild manifold $(M_\m, g_\m)$ with two ends. 

\vspace{0.2cm}

\item[(iii)] If $  \displaystyle  \lim_{ t \to 0 + } L(t)  = 0 $ and $ \m = 2 \c $, 
then  $(M, g)$ is isometric to the spatial Schwarzschild manifold $(M_\m, g_\m)$ with two ends. 

\end{enumerate} 

\end{theorem}

\begin{remark} \label{rem-ChengYau}
The metric $g$ is not assumed to be complete. 
If $g$ is complete with Ricci curvature bounded from below, 
then  $ \displaystyle \lim_{t \to 0} \max_{\Sigma_t} | \nabla u | = 0 $ 
by the gradient estimate of Cheng and Yau \cite{ChengYau75}.
This, together with the fact $  \int_{\Sigma_t} | \nabla u | = 4 \pi \c $, shows 
$ \displaystyle \lim_{ t \to 0 + }  L(t)   = 0$ under such assumptions. 
\end{remark}

\begin{remark}
For the purpose of the proof of the rigidity case in (ii), we record 
an equivalent form of the spatial Schwarzschild metric  $g_m$, $m \ne 0 $, 
written with respect to the level set gauge of the associated harmonic function.  
Upon a change of variable 
$$ t =  \left( 1 + \frac{ m}{2 \rho } \right)^{-1}, \ \text{where} \ \rho = | x | ,  $$
$g_m$ takes the form of  
\be \label{eq-gm-u-gauge}
g_m =   \left( 1 + \frac{ m}{2 | x |} \right)^4 g_0  
=  \frac{m^2}{4} \left[  \frac{ d t^2 }{   t^4 ( 1 - t)^4  } + \frac{ \sigma_o }{  t^2 ( 1 - t)^2 }  \right] ,
\ee
where $ \sigma_o$ denotes a round metric with Gauss curvature $1$ on the $2$-sphere $ S^2$. 
The sign of $m$ determines the domain of $t$: if $ m > 0 $, $ t \in (0,1)$; and if $ m < 0 $,
$ t \in (1, \infty)$. 

\vspace{.1cm}

Similarly, if $(M_r, g_m)$ is a portion of a spatial Schwarzschild manifold given in the introduction, 
via a further substitution of 
$$ s =  ( 1 - T)^{-1} ( t - T) , \ \ \text{where} \ T = \left( 1 + \frac{m}{2r} \right)^{-1} , $$ 
$g_m$ can be written as  
\be \label{eq-gm-u-gauge-r}
g_m = \frac{m^2}{4  ( 1 - T)^{2} }    \left[  \frac{ d s^2 }{  [( 1 - T) s + T]^4 (1 - s)^4  } 
  + \frac{\sigma_o }{  [( 1 - T) s + T]^2  (1 - s)^2 }    \right] .
\ee
This in particular implies 
\be \label{eq-L-t-u-gauge}
 \int_{s \, = \, \text{a constant} } | \nabla s |^2 \, d \sigma = 4 \pi 
    [( 1 - T) s + T]^{2} (1 - s)^{2}  , 
\ee
which gives \eqref{eq-L-t-S}.
\end{remark}

\begin{prof}[Proof of Theorem \ref{thm-app-1}] 
(i)  The monotone property of $\Q(t)$ can be obtained in a way similar 
to the derivation of \eqref{eq-Q-bd} in Remark \ref{rem-Q}. Below we give a different derivation 
that is better suited for showing the equality case. 

Given $(M, g)$ and $u$, consider the metric $ \bar g  = u^4 g $ on $M$, 
which may be incomplete, but 
has nonnegative scalar curvature and is asymptotically flat as $ | x |  \to \infty$. 

Given any $ T \in (0,1)$, a regular value of $ u $, 
let $ \Sigma_{_T} = u^{-1} (T)$. By the maximum principle, $ \Sigma_{_T} $ is a closed, connected surface, 
homologous to any coordinate sphere $ \{ |x| = \text{a constant} \}$. 

Let  $M_{_T} = \{ u \ge T \}$, which has boundary $ \Sigma_{_T}$. 
On $(M_{_T}, g )$, 
let $   u_{_T} = \frac{ u - T}{ 1 - T} $, 
which is harmonic, with $ u_{_T} \to 1 $ as $ |x| \to \infty$ and $ u_{_T} = 0 $ on $ \Sigma_{_T}$. 

On $(M_{_T}, \bar g)$, the function 
\be \label{eq-baru-T-T}
  \bar u_{_T}  = u^{-1} u_{_T}  = ( 1 - T)^{-1}  \left( 1 - T u^{-1} \right)   
\ee
is harmonic, with $ \bar u_{_T} \to 1 $ as $ |x| \to \infty$ and  $\bar u_{_T} = 0$ at $ \Sigma_{_T} $. 

Given any regular value $s \in [0,1)$ of $ \bar u_{_T}$, let 
 $$ \mathcal{ \bar B}_{_T}( s ) = \frac{1}{1 - s } \left[ 4 \pi -   \frac{1}{ (1 - s  )^2 } \int_{ \{ \bar u_{_T} = s \} } 
  | \bar \nabla \bar u_{_T} |_{\bar g} ^2 \, d \bar \sigma \right] ,
$$
where ``$ \bar \nabla $, $ | \cdot |_{\bar g}$, $ d \bar \sigma$" are the gradient, the length, the surface measure with respect to $\bar g$, respectively. 
By  \cite[Theorem 3.2 (ii)]{M22}, $ \mathcal{ \bar B}_{_T} ( s)$ is monotone nondecreasing in $s$, i.e. 
$$  \mathcal{ \bar B}_{_T}  ( s )  \,  \nearrow  \ \  \text{as} \  \ s \nearrow . $$
Given $s$, let $ t  \in [T,1) $ be chosen so that 
$ \{ \bar u_{_T} = s \} = \{ u  = t \} $.
By \eqref{eq-baru-T-T},
\be \label{eq-s-and-t}
s = \frac{1}{1 - T} \left( 1 - T t^{-1} \right) , \ \ 
1 - s  = \frac{T}{1- T} \frac{ 1 - t }{t} ,
\ee
and
\be \label{eq-conformal-bdry-u-ubar}
\int_{ \{ \bar u_{_T}  = s \} } | \bar \nabla \bar u_{_T} |_{\bar g}^2 \, d \bar \sigma = \int_{ \{ u = t \} } 
 | \nabla \bar u_{_T} |^2 \, d \sigma, \ \ \  
| \nabla \bar u_{_T}  | =   \frac{T }{1 - T} \, u^{-2} \, | \nabla u | . 
\ee
It follows from \eqref{eq-s-and-t} and  \eqref{eq-conformal-bdry-u-ubar} that 
\be \label{eq-bar-B-T-T}
\mathcal{\bar B}_{_T}  (s) = \mathcal{B}_{_T} (t), 
\ee
where 
\be \label{eq-bar-B-T-1-1}
\begin{split}
 \mathcal{B}_{_T}  (t)  
 = & \ \ \frac{1-T}{T} \,  \frac{t}{1 - t } \left[ 4 \pi -   \frac{1}{ t^2 (1 - t  )^2 } 
 \int_{ \{ u = t \} } \, | \nabla u |^2   \, d  \sigma \right]  \\
 = & \ \ \frac{1-T}{T} \, 4\pi  \, \Q (t) . 
 \end{split} 
\ee
As $ s = s (t)$ is strictly increasing in $t$, 
it follows from the monotonicity of $\mathcal{\bar B}_{_T} (s)$ that 
$$
\Q (t) \,   \nearrow  \  \text{as} \   t \nearrow , \ \forall \ t \in [T, 1). 
$$ 
As $ T \in (0, 1) $ is arbitrary, this shows the monotonicity of $\Q(t)$.

\vspace{.1cm}

By \cite[Theorem 3.2 (ii)]{M22}, $ \mathcal{\bar B}_{_T} (s )  $ also satisfies 
\be \label{eq-bar-u-T-T}
\mathcal{\bar B}_{_T} (s ) \le 4 \pi \m (\bar g) { \bar \c_{_T} }^{-1}  ,
\ee
where $ \m (\bar g)$ is the mass of $ \bar g$ and $\bar \c_{_T} $ is the coefficient in 
\bee
\bar u_{_T} = 1 - \bar \c_{_T} | x|^{-1} + o ( | x|^{-1} ) , \ \text{as} \ |x|  \to \infty. 
\eee
By the definitions of $\bar g$ and $\bar u_{_T}$, 
$ \m (\bar g) = \m - 2  \c $ and 
$  \bar \c_{_T} = \frac{T}{1-T} \, \c $. 
It follows from \eqref{eq-bar-B-T-T},  \eqref{eq-bar-B-T-1-1} and \eqref{eq-bar-u-T-T} that 
\be \label{eq-Bkt-0-0}
 \Q(t) \le \m \c^{-1} - 2 , \ \forall \ t \in [T, 1).
 \ee 
As $ T$ is arbitrary, \eqref{eq-Bkt-0-0} holds for any $ t \in (0,1)$.
 
Next, suppose \eqref{eq-Bkt-0-0} is an equality at some $ t \in (0,1)$. 
Choosing $ T = t $ in the above proof, we know  
  \eqref{eq-bar-u-T-T} becomes an equality at $ s = 0$. By  \cite[Theorem 3.2 (ii)]{M22}, 
this occurs if and only if  $(M_t , \bar g)$ is isometric to 
 $$ \left( \R^3 \setminus \{ | y |  < r \} , g_0 \right) $$
 for some $ r > 0$.
In terms of the coordinate $y$, functions $\bar u_t $ and $u$  on $M_t $ are given by 
\be \label{eq-bar-u-T-u-2}
  \bar u_t = 1 - \frac{   r}{ | y | } \ \  \text{and} \ \ 
  u^{-1} = 1 +  \frac{ ( 1 - t ) r}{ t | y | }    .  
\ee
Hence,  
$  g = u^{-4} g_0 = \left(  1 +  \frac{ (1-t) r}{  t  | y | }  \right)^4 g_0  . $
This shows $(M_t, g)$ is isometric to the part  
$$ \left(   \R^3 \setminus \{ | y |  < r \} , \left(  1 +  \frac{ \m  }{ 2   | y | } \right)^4 g_0 \right) $$
of a spatial Schwarzschild manifold of mass  
$ \m =  2  ( t^{-1} - 1 )  r   > 0   $.

Combined with \eqref{eq-c-m-S}, i.e. 
$  \c_{_{\Sigma_t}}  = r + \frac{\m}{2} $, 
and the fact 
$  \c_{_{\Sigma_t} } = \frac{1}{1-t} \c $, the above in particular shows $ \m = 2 \c $ in this case. 

\vh

(ii)  By the properties of $\Q (t)$ in (i), 
 the limit 
 $$ \lim_{t \to 0 +} \Q(t) = - \lim_{t \to 0 + } \frac{L(t)}{t} \ \ \text{exists} \  (\text{finite or $\infty$}) , $$
 and satisfies 
 \be \label{eq-Bkt-bar-0}
-   \lim_{t \to 0} \frac{L(t)}{t} \le \Q (t) \le   \m \c^{-1}  - 2  , \ \forall \, t . 
 \ee 
Consequently, the equality 
 $$ -   \lim_{t \to 0} \frac{L(t)}{t} =   \m \c^{-1}  - 2 $$ 
 holds 
 if and only if 
\be \label{eq-q-t-rigidity}
 \Q (t) = \m \c^{-1}  - 2, \  \forall \, t . 
\ee
By the rigidity conclusion in (i), this means $\m > 0 $ and, for every $ t \in (0,1)$, 
$(M_t, g)$ is isometric to 
\be \label{eq-part-M-m-with-rt-known}
 \left(   \R^3 \setminus \left\{ | y |  <  r_t \right\} , \left(  1 +  \frac{ \m  }{ 2   | y | } \right)^4 g_0 \right) , 
\ \text{with} \  r_t = \frac{ \m t}{2 (1-t)}   .  
\ee  
To show this  implies $(M,g)$ itself  is isometric to $(M_\m, g_\m)$, we make use of the 
level set gauge of $u$ in $(M, g)$.

\vspace{.1cm}

By the formula of $u$ in \eqref{eq-bar-u-T-u-2}, on $M_t = \{ u \ge t \}$, $u$ satisfies 
\be \label{eq-grad-u} 
 | \nabla u |  
 =  \frac{2}{ \m }   u^2 (1-u)^2   > 0 .
\ee 
As $ t $ is arbitrary, $ \nabla u $ never vanishes on $M$. As a result, 
by considering the integral flow of $ | \nabla u |^{-2}  \nabla u $, 
$(M,g)$ is isometric to 
\be \label{eq-level-set-gauge}
\begin{split}
 & \ \left( (0, 1) \times S^2 , \frac{1}{ | \nabla u |^2} d t^2 + g_t \right) \\
=  & \ \left( (0, 1) \times S^2  ,  \frac{\m^2}{4}  \frac{1}{t^4 ( 1 - t)^4} \, d t^2+ g_t  \right) . 
 \end{split}
\ee 
Here we have used the fact that, for $ t $ close to $ 1$, $\Sigma_t$ is a $2$-sphere $ S^2$. 
The metrics $ \{g_t\}_{0 < t < 1} $ is the family of metrics 
induced from $g$ on $\Sigma_t = \{ t \} \times S^2$.  
The derivative of $g_t$ with respect to $t$ satisfies 
\be \label{eq-1st-fundamental-form}
 \frac{d }{d t} g_t  = \frac{2}{| \nabla u |} \Pi_t =  \frac{ \m }{  t^2 (1-t)^2 } \Pi_t  , 
\ee
where $\Pi_t$ is the second fundamental form of $\Sigma_t $ in $(M,g)$ with respect to $  | \nabla u|^{-1} \nabla u$. 

For each $ t \in (0,1)$, by  \eqref{eq-gm-u-gauge} and \eqref{eq-part-M-m-with-rt-known}, 
$(M_t, g)$ is isometric to  
\be \label{eq-key-part-in-rigidity}
 ([t, 1) \times S^2, g_\m),
\ee
where
$$
g_\m = \frac{\m^2}{4}  \left[  \frac{1 }{   s^4 ( 1 - s )^4  }  d s^2 + \frac{1 }{  s^2 ( 1 - s)^2 }  \sigma_o  \right] .
$$
In particular, this implies $ g_t $ has constant Gauss curvature 
$ \displaystyle  \frac{4 } { \m^2} t^2 (1-t)^2 $, 
and the second fundamental form 
$ \Pi_t$ satisfies 
\be \label{eq-2nd-fundamental-form-1}
 \frac{d}{ds} |_{s = t}  \left(  \frac{ \m^2  }{ 4 s^2 ( 1 - s)^2 }  \sigma_o \right) 
 =  \frac{ \m  }{   t^2 ( 1 - t )^2  }   \Pi_t .
\ee 
By \eqref{eq-2nd-fundamental-form-1},  $\Pi_t$ is a function multiple of $g_t$;  
more precisely, 
\be \label{eq-2nd-fundamental-form-2}
 \frac{d}{d s} |_{s = t}  \ln  \left(   s^{-2}  ( 1 - s)^{-2}   \right)  \, g_t 
 =  \frac{ \m  }{   t^2 ( 1 - t )^2  }   \Pi_t .
\ee 
It follows from \eqref{eq-1st-fundamental-form} and \eqref{eq-2nd-fundamental-form-2} that
\be \label{eq-ode}
 \frac{d }{d t} g_t   =    \frac{d}{dt} \ln \left[   t^{-2} ( 1 - t)^{-2}   \right] \, g_t . 
\ee
At $ t = \frac12$, the metric $g_\frac12$ has constant Gauss curvature  $ ( 4 \m^2)^{-1}$, 
and hence
\be \label{eq-initial-condition}
 g_{\frac12} = 4 \m^2 \tau_o ,
\ee
where $ \tau_o$ is a round metric of Gauss curvature $1$ on $ S^2$. 
It follows from \eqref{eq-ode} and \eqref{eq-initial-condition} that 
$$ g_t =  \frac{\m^2}{4 t^2 ( 1 - t)^2 } \, \tau_o, \ \forall \, t \in (0,1) . $$
This, together with \eqref{eq-level-set-gauge} and \eqref{eq-gm-u-gauge}, shows $(M,g)$ is isometric to $(M_\m, g_\m)$.

\vh 

(iii) As before, given any  regular value $ T \in (0,1) $ of $ u $, 
let $ \Sigma_{_T} = u^{-1} (T)$ and   $M_{_T} = \{ u \ge T \}$.
On $(M_{_T}, g )$, let $ u_{_T} = \frac{ u - T}{ 1 - T} $. 
Let $ \c_{_T}$ denote the capacity of $ \Sigma_{_T}$ in $(M_{_T} ,g)$, then
$ \c_{_T} = ( 1 - T)^{-1} \c $. 

Suppose $ \m \le 2 \c $, then 
$$ \m  < 2 \c_{_T}  . $$   
Theorem \ref{thm-main-2} is hence applicable to the tuple  $(M_{_T}, g, u_{_T})$.
Therefore, for any regular value $ s \in [0,1)$ of $u_{_T}$, 
the function 
$$ \mathcal{S}_{_T} (s) = 2 \left[ \frac{s}{1-s}  
- \frac{1}{(1-s)^2}  \left( \frac{1}{4\pi} \int_{ \{ u_{_T} = s \} } | \nabla u_{_T} |^2 \, d \sigma \right)^\frac12 \right] $$
is monotone nondecreasing in $s$ and 
\be \label{eq-S-T-m-c}
  \mathcal{S}_{_T} (s) \le \m \c_{_T}^{-1} - 2 . 
\ee
Let $ t = ( 1 - T) s + T$, then 
$$ \mathcal{S}_{_T} (s) + 2 = ( 1 - T ) \left(  \mathcal{S}(t) + 2  \right) ,  $$
where
$$ \mathcal{S}(t) = 
2 \left[ \frac{t }{1-t}  
- \frac{ 1 }{(1-t)^2}   \left( \frac{1}{4\pi} \int_{ \{ u  = t  \} } | \nabla u |^2 \, d \sigma \right)^\frac12 \right]. 
$$
Thus, by the properties of $\mathcal{S}_{_T} (s)$, 
$  \mathcal{S}(t) $ is monotone nondecreasing in $t \in [ T, 1) $,  and 
\be \label{eq-S-t-m-c}
\mathcal{S}(t) \le \m \c^{-1} - 2 . 
\ee 
As $ T \in (0,1) $ is arbitrary, the above holds for any $t  \in (0,1)$.

\vspace{.1cm}

Letting $ t \to 0 +$ in \eqref{eq-S-t-m-c} gives 
\be \label{eq-sandwich} 
 -2 \lim_{t \to 0 +} L(t)^\frac12 
\le  \mathcal{S}(t) \le \m \c^{-1} - 2, \ \forall \, t \in (0,1). 
\ee
If  $  \displaystyle \lim_{ t \to 0 + } \int_{ \Sigma_t }  | \nabla u |^2 \, d \sigma  = 0$ and $ \m = 2 \c $, 
then, by \eqref{eq-sandwich},  
$$   \mathcal{S}(t)  = \m \c^{-1} - 2 = 0 , \ \forall \, t \in (0,1).  $$
Consequently, 
$$ L (t) = t^2 (1 - t)^2 , \  \forall \, t \in (0,1).$$ 
In particular, 
$$ \lim_{t \to 0 +} \frac{L(t) }{t} = 0 . $$
It follows from the rigidity conclusion in (ii) that $(M, g)$ is isometric to $(M_\m, g_\m)$. 
\end{prof}






%
%


\begin{thebibliography}{99}


\bibitem{AMO21}
Agostiniani, V., 
Mazzieri, L.,  
Oronzio, F.:
A Green’s function proof of the positive mass theorem. 
\emph{Commun. Math. Phys.}, \textbf{405}, no. 2, Paper No. 54, 23 pp, (2024).

\bibitem{AMMO22}
Agostiniani, V., 
Mantegazza, C., 
Mazzieri, L., 
Oronzio, F.: 
Riemannian Penrose inequality via nonlinear potential theory. 
\emph{Rend. Lincei Mat. Appl.},  \textbf{34},  715--726 (2023)

\bibitem{ADM61}
Arnowitt, R., 
Deser, S., 
Misner, C. W.:
Coordinate invariance and energy expressions in general relativity. 
\emph{Phys. Rev.}, \textbf{122}, no. 3, 997--1006 (1961)

\bibitem{Bartnik86} 
Bartnik, R.:  
The mass of an asymptotically flat manifold. 
\emph{Comm. Pure Appl. Math.},  \textbf{39}, no. 5, 661--693  (1986) 

\bibitem{BKKS19}
Bray,  H.L.,  
Kazaras, D.P., 
Khuri, M.A., 
Stern, D.L.: 
Harmonic functions and the mass of $3$-dimensional asymptotically flat Riemannian manifolds. 
\emph{J. Geom. Anal.},  \textbf{32}, no. 6, Paper No.184, 29 pp (2022)

\bibitem{ChengYau75}
Cheng, S.-Y, 
Yau, S.-T.:  
Differential equations on Riemannian manifolds and their geometric applications.  
\emph{Comm. Pure Appl. Math.},  \textbf{28}, no. 3, 333–354 (1975)

\bibitem{CL21}
Chodosh, O., 
Li, C.:
Stable minimal hypersurfaces in $\R^4$. 
\emph{Acta Math.}, \textbf{233}, no. 1, 1--31 (2024)

\bibitem{Chrusciel86}
Chru\'sciel, P.: 
Boundary conditions at spatial infinity from a Hamiltonian point of view.
\emph{Topological Properties and Global Structure of Space-Time}, Plenum Press, New York, 49--59 (1986)

\bibitem{Colding12}
Colding, T.H.:  
New monotonicity formulas for Ricci curvature and applications I.
\emph{Acta Math.},  \textbf{209}, no. 2, 229--263 (2012)

\bibitem{CM13}
Colding, T.H.,   
Minicozzi,  W.P.: 
Monotonicity and its analytic and geometric implications. 
\emph{Proc. Natl. Acad. Sci.},  USA \textbf{110}, no. 48, 19233--19236 (2013)


\bibitem{CM14}
Colding, T.H.,   
Minicozzi,  W.P.:  
Ricci curvature and monotonicity for harmonic functions. 
\emph{Calc. Var. Partial Differential Equations},  \textbf{49}, no. 3-4, 1045--1059  (2014)

\bibitem{M22}
Miao, P.:
Mass, capacitary functions, and the mass-to-capacity ratio. 
\emph{Peking Math J.}, \textbf{8}, 351--404  (2025)

\bibitem{M24}
Miao, P.: 
Implications of some mass-capacity inequalities. 
\emph{J. Geom. Anal.}, \textbf{34},  Paper No. 241, 16 pp (2024)

\bibitem{MW21}
Munteanu, O., 
Wang, J.:
Comparison theorems for three-dimensional manifolds with scalar curvature bound. 
\emph{Int. Math. Res. Not.}, \textbf{2023}, no. 3, 2215--2242 (2023)

\bibitem{MW22}
Munteanu, O., 
Wang, J.:
Geometry of three-dimensional manifolds with positive scalar curvature.
\emph{Amer. J. Math.}, to appear, arXiv:2201.05595v3 (2022)

\bibitem{O22}
Oronzio, F.: 
ADM mass, area and capacity in asymptotically flat $3$-manifolds with nonnegative scalar curvature. 
\emph{Commun. Contemp. Math.}, \textbf{27},  no. 9, 2550011 (2025)

\bibitem{SchoenYau79} 
Schoen, R., 
Yau,  S.-T.:  
On the proof of the positive mass conjecture in general relativity.
\emph{Commun.  Math. Phys.},  \textbf{65}, no. 1, 45--76 (1979)

\bibitem{Stern19}
Stern, D.: 
Scalar curvature and harmonic maps to ${S}^1$.  
\emph{J. Differ. Geom.}, \textbf{122}, no. 2, 259--269 (2022)

\bibitem{XYZ23}
Xia, C., 
Yin, J., 
Zhou, X.:
New monotonicity for $p$-capacitary functions in $3$-manifolds with nonnegative scalar curvature. 
\emph{Adv. Math.}, \textbf{440}, 109526 (2024)

\bibitem{Witten81} 
Witten, E.: 
A new proof of the positive energy theorem.
\emph{Commun.  Math. Phys.}, \textbf{80}, no. 3, 381--402  (1981)

\end{thebibliography}
\end{document}